# THE α-CALCULUS-CUM-α-ANALYSIS OF $\dfrac{\partial^r}{\partial s^r}\zeta(s,\alpha)$

V.V. RANE

A-3/203, ANAND NAGAR,

DAHISAR, MUMBAI-400 068,

INDIA

v_v_rane@yahoo.co.in

**Abstract** : For Hurwitz zeta function $\zeta(s,\alpha)$, as a function of $\alpha$, we discuss the location and the nature of singularities of $\zeta^{(r)}(s,\alpha) = \dfrac{\partial^r}{\partial s^r}\zeta(s,\alpha)$, the formulae for the derivatives and the primitives of $\zeta^{(r)}(s,\alpha)$, the Riemann-integrability of $\zeta^{(r)}(s,\alpha)$ on the intervals [0,1] and $[1,\infty)$; and the evaluation of integrals
$$\int_0^1 \zeta(-m_1,\alpha)\cdot\zeta(-m_2,\alpha)............\zeta(-m_k,\alpha)\cdot\zeta^{(r)}(s,\alpha)d\alpha\ ,\ \text{where}\ m_1,m_2,............,m_k \geq 0$$
are integers and Re s<1.

# THE α-CALCULUS-CUM-α-ANALYSIS OF $\frac{\partial^r}{\partial s^r}\zeta(s,\alpha)$

V.V.RANE

A-3/203, ANAND NAGAR,

DAHISAR, MUMBAI-400 068,

INDIA

v_v_rane@yahoo.co.in

Let $r \geq 0$ be an integer and for the complex variables $s, \alpha$, let $\zeta(s,\alpha)$ be the Hurwitz zeta function and let $\zeta^{(r)}(s,\alpha) = \frac{\partial^r}{\partial s^r}\zeta(s,\alpha)$. The object of this paper is to bring out the behaviour of $\zeta^{(r)}(s,\alpha)$ as an analytic function of the complex variable $\alpha$; to determine the location and the nature of the singularities of $\zeta^{(r)}(s,\alpha)$; to determine the Riemann integrability of $\zeta^{(r)}(s,\alpha)$ on intervals $[0,1]$ and $[1,\infty)$; to obtain the formulae for its derivatives and primitives with respect to $\alpha$; to evaluate the integrals $\int_0^1 \zeta(-m_1,\alpha)\zeta(-m_2,\alpha)\ldots\ldots\ldots\zeta(-m_r,\alpha)\zeta(s,\alpha)$ and $\int_0^1 \zeta(-m_1,\alpha)\cdot\zeta(-m_2,\alpha)\ldots\ldots\ldots\zeta(-m,\alpha)\cdot\zeta'(s,\alpha)d\alpha$, where $m_1, m_2, \ldots\ldots\ldots, m_r \geq 0$ are integers and Re s<1.

Next, we formally introduce our notation and terminology. For complex $\alpha \neq 0, -1, -2, \ldots\ldots\ldots$, and for the complex variable s, let $\zeta(s,\alpha)$ be the Hurwitz zeta function defined by $\zeta(s,\alpha) = \sum_{n \geq 0}(n+\alpha)^{-s}$ for Re s>1 and its analytic continuation. Let $\zeta(s,1) = \zeta(s)$, the Riemann zeta function. In what follows, $\Gamma(s)$ stands for gamma function. $B_n(\alpha) = \sum_{i=0}^{n} \binom{n}{i} B_{n-i} \alpha^i$ stands for Bernoulli polynomial of degree n. Here $B_n = B_n(0)$ are Bernoulli numbers, which are known to be rational numbers. We note that if $n \geq 0$ is an integer, then $\zeta(-n,\alpha) = -\frac{B_{n+1}(\alpha)}{n+1}$ and $B_n(1-\alpha) = (-1)^n B_n(\alpha)$. We write $\psi(\alpha) = \frac{\Gamma'}{\Gamma}(\alpha)$. Note that $\frac{d}{d\alpha}\zeta'(0,\alpha) = \frac{d}{d\alpha}\log\frac{\Gamma(\alpha)}{\sqrt{2\pi}} = \frac{\Gamma'}{\Gamma}(\alpha) = \psi(\alpha)$. Thus $\frac{d}{d\alpha}\zeta'(0,\alpha) = \psi(\alpha)$. We shall see that $\frac{d}{d\alpha}\psi(\alpha) = \zeta(2,\alpha)$. In view of the fact $\frac{d}{d\alpha}\zeta(s,\alpha) = -s\zeta(s+1,\alpha)$, we have $\frac{d^r}{d\alpha^r}\psi(\alpha) = (-1)^{r-1} r! \zeta(r+1,\alpha)$. We also note that if $s_1, s_2$ are complex numbers, then

$$\int_0^1 \zeta(s_1,\alpha) \cdot \zeta(s_2,\alpha) d\alpha = 2(2\pi)^{s_1+s_2-2} \cdot \Gamma(1-s_1)\Gamma(1-s_2) \cdot \cos\frac{\pi}{2}(s_1-s_2) \cdot \zeta(2-s_1-s_2),$$

barring the singularities of either side. In author [3], we study evaluation of Tornheim double zeta function $T(s_1,s_2,s_3)$ for various values of the complex variables $s_1,s_2,s_3$ in terms of integrals of the type $\int_0^1 \zeta(s_1,\alpha)\zeta(s_2,\alpha)\zeta^{(r)}(s_3,\alpha)d\alpha$ for $r = 0$ or 1. In view of this, we shall find some integrals of the type $\int_0^1 \zeta(s_1,\alpha)\zeta(s_2,\alpha)\zeta^{(r)}(s_3,\alpha)d\alpha$, which are explicitly computable in terms of values of Riemann zeta function and its derivatives. Apart from several other facts, in particular, we prove the following.

**Proposition :** For any integer $r \geq 0$ and for Re s<o, $\zeta^{(r)}(s,\alpha)$ is a continuous function of $\alpha$ in the whole complex $\alpha$-plane and we have the following.

1 ) $\zeta^{(r)}(s,0) = \zeta^{(r)}(s) = \zeta^{(r)}(s,1)$ for Re s<o and $\lim_{s\to 0-} \zeta^{(r)}(s,0) = \zeta^{(r)}(0)$,

where $s \to 0$ through real values from the left of $0$.

2) $\dfrac{\partial}{\partial \alpha} \zeta^{(r)}(s,\alpha) = -r \cdot \zeta^{(r-1)}(s+1,\alpha) - s\zeta^{(r)}(s+1,\alpha)$ for $s \neq 0$ and for $r \geq 0$.

3) We have $\dfrac{\partial}{\partial \alpha} \zeta^{(r)}(0,\alpha) = -r! \gamma_{r-1}(\alpha)$ for $r \geq 0$,

where $s\zeta(s+1,\alpha) = \sum_{n \geq 0} \gamma_{n-1}(\alpha) s^n$ with $\gamma_{-1}(\alpha) = 1$.

4) We have $\dfrac{\partial}{\partial \alpha} \gamma_{r-1}(\alpha) = -\dfrac{1}{r!} \dfrac{\partial^2}{\partial \alpha^2} \zeta^{(r)}(0,\alpha) = -\dfrac{1}{r!} \cdot \dfrac{\partial^r}{\partial s^r} \left( s(s+1)\zeta(s+2,\alpha) \right) \Big|_{s=0}$ for $r \geq 1$.

**Note : 1) of** Proposition follows from the fact $\zeta(s,\alpha) - \zeta(s,\alpha+1) = \alpha^{-s}$

and consequently, $\zeta^{(r)}(s,\alpha) - \zeta^{(r)}(s,\alpha+1) = \alpha^{-s}(-\log \alpha)^r$ for an integer $r \geq 0$,

after letting $\alpha \to 0$.

**Corollaries of Proposition :**

1) For $s \neq 1$, the primitives (or the indefinite integrals) of $\zeta(s,\alpha), \zeta'(s,\alpha)$ and

$\zeta''(s,\alpha)$ are as follows :

a) $\int \zeta(s,\alpha) d\alpha = \dfrac{\zeta(s-1,\alpha)}{1-s}$

b) $\int \zeta'(s,\alpha) d\alpha = \dfrac{\zeta(s-1,\alpha)}{(1-s)^2} + \dfrac{\zeta'(s-1,\alpha)}{1-s}$

c) $\int \zeta''(s,\alpha) d\alpha = 2\dfrac{\zeta(s-1,\alpha)}{(1-s)^3} + 2\dfrac{\zeta'(s-1,\alpha)}{(1-s)^2} + \dfrac{\zeta''(s-1,\alpha)}{1-s}$

d) $\int \zeta(s,1-\alpha) d\alpha = \dfrac{\zeta(s-1,1-\alpha)}{s-1}$

e) $\int \zeta'(s,1-\alpha) d\alpha = -\left( \dfrac{\zeta(s-1,1-\alpha)}{(1-s)^2} + \dfrac{\zeta'(s-1,1-\alpha)}{1-s} \right)$

f) In general, for an integer $r \geq 0$,

$$\int \zeta^{(r)}(s,\alpha)d\alpha = \sum_{\ell=0}^{r} c_\ell \frac{\zeta^{(\ell)}(s-1,\alpha)}{(1-s)^{r+1-\ell}}$$

for some absolute constants $c_\ell$'s.

2) Under the action of the operator $\frac{\partial}{\partial \alpha}$, we have the following diagram

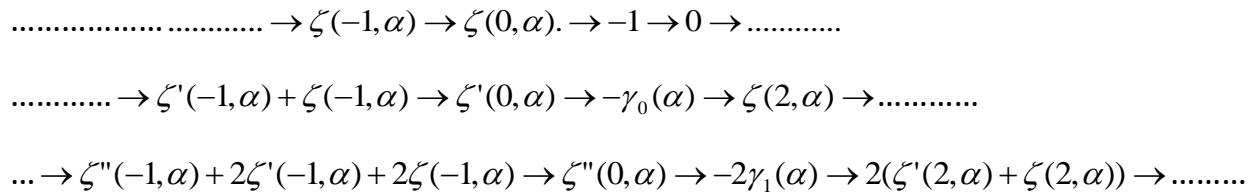

$$\ldots \to \zeta(-1,\alpha) \to \zeta(0,\alpha) \to -1 \to 0 \to \ldots$$

$$\ldots \to \zeta'(-1,\alpha) + \zeta(-1,\alpha) \to \zeta'(0,\alpha) \to -\gamma_0(\alpha) \to \zeta(2,\alpha) \to \ldots$$

$$\ldots \to \zeta''(-1,\alpha) + 2\zeta'(-1,\alpha) + 2\zeta(-1,\alpha) \to \zeta''(0,\alpha) \to -2\gamma_1(\alpha) \to 2(\zeta'(2,\alpha) + \zeta(2,\alpha)) \to \ldots$$

3) We have $\int_1^\infty \zeta^{(r)}(s,\alpha)d\alpha = -\sum_{\ell=o}^{r} c_\ell \cdot \frac{\zeta^{(\ell)}(s-1)}{(1-s)^{r+1-\ell}}$ for Re s>2,

where $c_\ell$'s are absolute constants as in f) of Corollary 1) above.

**Note**: This is so, because for Re s>1, as $\alpha \to \infty$,

$$\zeta^{(r)}(s,\alpha) = (-1)^r \sum_{n \geq 0}(n+\alpha)^{-s} \log^r(n+\alpha) = 0$$

4) We have $\int_0^1 \zeta^{(r)}(s,\alpha)d\alpha = 0$ for Re s<1.

**Note**: This follows from f) of corollary 1) and 1) of Proposition above.

5) For Re $s_1$, Re $s_2$ <0 and for real s, we have

$$\lim_{s \to 1-} \int_0^1 \zeta(s_1,\alpha) \cdot \zeta(s_2,\alpha) \cdot (s-1)\zeta(s,\alpha)d\alpha = \int_0^1 \zeta(s_1,\alpha) \cdot \zeta(s_2,\alpha)d\alpha - \zeta(s_1) \cdot \zeta(s_2)$$

$$= 2(2\pi)^{s_1+s_2-2} \cdot \Gamma(1-s_1)\Gamma(1-s_2)\cos\frac{\pi}{2}(s_1-s_2) \cdot \zeta(2-s_1-s_2) - \zeta(s_1) \cdot \zeta(s_2),$$

where $s \to 1$ from left through real values.

6) Let $m_1, m_2, \ldots, m_r \geq 1$ be integers and let $N = \sum_{i=1}^{r} m_i$. Then

$$\int_0^1 \zeta(1-m_1, \alpha) \cdot \zeta(1-m_2, \alpha) \ldots \zeta(1-m_r, \alpha) d\alpha$$

$$= \frac{(-1)^r}{m_1 m_2 \ldots m_r} \int_0^1 B_{m_1}(\alpha) \cdot B_{m_2}(\alpha) \ldots B_{m_r}(\alpha) d\alpha \text{ is explicitly computable as}$$

a rational number, which equals zero when N is odd.

**Remark :** Note that $\prod_{\ell=1}^{r} B_{m_\ell}(\alpha)$ is a polynomial with rational coefficients.

Hence $\int_0^1 \prod_{\ell=1}^{r} B_{m_\ell}(\alpha) d\alpha$ is a rational number.

Note that $\int_0^1 B_{m_1}(\alpha) B_{m_2}(\alpha) \ldots B_{m_r}(\alpha) d\alpha$

$$= \int_0^1 B_{m_1}(1-\alpha) \cdot B_{m_2}(1-\alpha) \ldots B_{m_r}(1-\alpha) d\alpha = (-1)^N \int_0^1 B_{m_1}(\alpha) \cdot B_{m_2}(\alpha) \ldots B_{m_r}(\alpha) d\alpha$$

=0 , when N is odd.

7) If $m_1, m_2, \ldots, m_r \geq 0$ are integers and Re s<1, then

$$\int_0^1 \zeta(-m_1, \alpha) \cdot \zeta(-m_2, \alpha) \ldots \cdot \zeta(-m_r, \alpha) \zeta(s, \alpha) d\alpha$$

is explicitly computable as a linear combination of

$\zeta(s-1), \zeta(s-2), \ldots, \zeta(s-N)$ with coefficients dependent on s,

where $N = \sum_{i=1}^{r}(m_i + 1)$ is the degree of the product polynomial

$\zeta(-m_1, \alpha) \cdot \zeta(-m_2, \alpha) \ldots \cdot \zeta(-m_r, \alpha)$.

8) If $m_1, m_2, \ldots, m_r \geq 0$ are integers and Re s<1,

then $\int_0^1 \zeta(-m_1,\alpha)\cdot\zeta(-m_2,\alpha)............\zeta(-m_r,\alpha)\cdot\zeta'(s,\alpha)d\alpha$

is explicitly computable as a linear combination of

$\zeta'(s-1),\zeta'(s-2),..............,\zeta'(s-N)$ ; $\zeta(s-1),\zeta(s-2),..............,\zeta(s-N)$

with coefficients dependent upon s , where $N = \sum_{i=1}^{r}(m_i +1)$.

9) We have for Re s>1 ,

$$\int_0^1 \zeta(0,\alpha)\zeta(1-s,\alpha)\cdot\zeta(2-s,\alpha)d\alpha = \frac{1}{2(s-1)}\left(2(2\pi)^{-2s}\cdot\Gamma^2(s)\zeta(2s)-\zeta^2(1-s)\right)$$

      Assuming Proposition , we shall prove the corollaries .

**Proof of Corollary 2)** : We shall sketch the proof of

$$\frac{\partial}{\partial\alpha}(-\gamma_0(\alpha)) = \frac{\partial^2}{\partial\alpha^2}\zeta'(0,\alpha) = \zeta(2,\alpha).$$

We have $\frac{\partial^2}{\partial\alpha^2}\zeta'(0,\alpha) = \frac{\partial^2}{\partial\alpha^2}\frac{\partial}{\partial s}\zeta(s,\alpha)|_{s=0}$

$= \frac{\partial}{\partial s}\frac{\partial^2}{\partial\alpha^2}\zeta(s,\alpha)|_{s=0} = \frac{\partial}{\partial s}\frac{\partial}{\partial\alpha}(-s\zeta(s+1,\alpha))|_{s=0} = \frac{\partial}{\partial s}(s(s+1)\zeta(s+2,\alpha))|_{s=0}$

$= ((s+1)\zeta(s+2,\alpha) + s\zeta(s+2,\alpha) + s(s+1)\zeta'(s+2,\alpha))_{s=0} = \zeta(2,\alpha)$ .

      Next , we show $\frac{\partial}{\partial\alpha}\zeta''(0,\alpha) = -2\gamma_1(\alpha)$

and $\frac{\partial}{\partial\alpha}(-2\gamma_1(\alpha)) = \frac{\partial^2}{\partial\alpha^2}\zeta''(0,\alpha) = 2(\zeta'(2,\alpha)+\zeta(2,\alpha))$ .

We have $\frac{\partial}{\partial\alpha}\zeta''(0,\alpha) = \frac{\partial}{\partial\alpha}\frac{\partial^2}{\partial s^2}\zeta(s,\alpha)|_{s=0}$

$= \frac{\partial^2}{\partial s^2}\frac{\partial}{\partial\alpha}\zeta(s,\alpha)|_{s=0} = \frac{\partial^2}{\partial s^2}(-s\zeta(s+1,\alpha))|_{s=0} = -(s\zeta(s+1,\alpha))''|_{s=0} = -2\gamma_1(\alpha)$ .

Next, $\dfrac{\partial}{\partial \alpha}(-2\gamma_1(\alpha)) = \dfrac{\partial^2}{\partial \alpha^2}\zeta''(0,\alpha) = \dfrac{\partial^2}{\partial \alpha^2}\dfrac{\partial^2}{\partial s^2}\zeta(s,\alpha)\big|_{s=0}$

$= \dfrac{\partial^2}{\partial s^2}\dfrac{\partial^2}{\partial \alpha^2}\zeta(s,\alpha)\big|_{s=0} = \dfrac{\partial^2}{\partial s^2}(s(s+1)\zeta(s+2,\alpha))\big|_{s=0}$

$= \dfrac{\partial}{\partial s}\{s\cdot\zeta(s+2,\alpha) + (s+1)\zeta(s+2,\alpha) + s(s+1)\zeta'(s+2,\alpha)\}_{s=0}$

$= \begin{cases}\zeta(s+2,\alpha) + s\zeta'(s+2,\alpha) + \zeta(s+2,\alpha) + (s+1)\zeta'(s+2,\alpha) \\ + (s+1)\zeta'(s+2,\alpha) + s\zeta'(s+2,\alpha) + s(s+1)\zeta''(s+2,\alpha)\end{cases}\bigg|_{s=0} = 2(\zeta(2,\alpha) + \zeta'(2,\alpha))$.

**Proof of Corollary 5 ) :** For Re $s_1$, Re $s_2$ <0 and for real s with s<1 ,

$\displaystyle\int_0^1 \zeta(s_1,\alpha)\cdot\zeta(s_2,\alpha)\cdot(s-1)\zeta(s,\alpha)d\alpha = -\int_0^1 \zeta(s_1,\alpha)\cdot\zeta(s_2,\alpha)\dfrac{\partial}{\partial\alpha}\zeta(s-1,\alpha)d\alpha$

$= -\left\{\zeta(s_1,\alpha)\cdot\zeta(s_2,\alpha)\cdot\zeta(s-1,\alpha)\big|_{\alpha=0}^{1} - \displaystyle\int_0^1 \zeta(s-1,\alpha)\cdot\dfrac{\partial}{\partial\alpha}(\zeta(s_1,\alpha)\cdot\zeta(s_2,\alpha))d\alpha\right\}$

$= 0 + \displaystyle\int_0^1 \zeta(s-1,\alpha)\cdot\dfrac{\partial}{\partial\alpha}(\zeta(s_1,\alpha)\cdot\zeta(s_2,\alpha))d\alpha$ .

Thus $\displaystyle\lim_{s\to 1-}\int_0^1 \zeta(s_1,\alpha)\cdot\zeta(s_2,\alpha)(s-1)\zeta(s,\alpha)d\alpha$

$= \displaystyle\lim_{s\to 1-}\int_0^1 \zeta(s-1,\alpha)\cdot\dfrac{\partial}{\partial\alpha}(\zeta(s_1,\alpha)\cdot\zeta(s_2,\alpha))d\alpha = \int_0^1 \zeta(0,\alpha)\cdot\dfrac{\partial}{\partial\alpha}(\zeta(s_1,\alpha)\cdot\zeta(s_2,\alpha))d\alpha$

$= \displaystyle\int_0^1 (\dfrac{1}{2}-\alpha)\dfrac{\partial}{\partial\alpha}(\zeta(s_1,\alpha)\cdot\zeta(s_2,\alpha))d\alpha = (\dfrac{1}{2}-\alpha)\zeta(s_1,\alpha)\cdot\zeta(s_2,\alpha)\big|_{\alpha=0}^1 + \int_0^1 \zeta(s_1,\alpha)\cdot\zeta(s_2,\alpha)d\alpha$ ,

on integration by parts .

This , in turn , $= \displaystyle\int_0^1 \zeta(s_1,\alpha)\cdot\zeta(s_2,\alpha)d\alpha - \zeta(s_1)\cdot\zeta(s_2)$

$= 2(2\pi)^{s_1+s_2-2}\cdot\Gamma(1-s_1)\Gamma(1-s_2)\cos\dfrac{\pi}{2}(s_1-s_2)\cdot\zeta(2-s_1-s_2) - \zeta(s_1)\cdot\zeta(s_2)$.

**Proof of Corollary 7) :** We have

$$\zeta(-m_1,\alpha)\cdot\zeta(-m_2,\alpha)\ldots\ldots\ldots\zeta(-m_r,\alpha) = \left(-\frac{B_{m_1+1}(\alpha)}{m_1+1}\right)\left(-\frac{B_{m_2+1}(\alpha)}{m_2+1}\right)\ldots\ldots\ldots\left(-\frac{B_{m_r+1}(\alpha)}{m_r+1}\right)$$

$= P_N(\alpha)$, say, where $P_N(\alpha)$ is a polynomial of degree $N = \sum_{i=1}^{r}(m_i+1)$.

Let $P_N(\alpha) = \sum_{\ell=0}^{N} a_i \alpha^i$, where $a_i$'s are rational numbers.

Thus $\int_0^1 P_N(\alpha)\zeta(s,\alpha)d\alpha = \sum_{i=0}^{N} a_i \int_0^1 \alpha^i \zeta(s,\alpha)d\alpha$.

Note that for each i , $\int_0^1 \alpha^i \zeta(s,\alpha) d\alpha = \int_0^1 \alpha^i d\left(\frac{\zeta(s-1,\alpha)}{1-s}\right)$

$= \alpha^i \frac{\zeta(s-1,\alpha)}{1-s}\Big|_{\alpha=0}^{1} - \int_0^1 \frac{\zeta(s-1,\alpha)}{1-s} i\alpha^{i-1} d\alpha = \frac{\zeta(s-1)}{1-s} + \frac{i}{s-1}\int_0^1 \alpha^{i-1}\zeta(s-1,\alpha)d\alpha$ .

We write $\int_0^1 \alpha^{i-1} \zeta(s-1,\alpha) d\alpha = \int_0^1 \alpha^{i-1} d\left(\frac{\zeta(s-2,\alpha)}{2-s}\right)$ and continue the process of integration by parts till we reach the stage , when we come across the integral

$\int_0^1 \zeta(s-i,\alpha) d\alpha$, which is equal to zero .

Thus $\int_0^1 \alpha^i \zeta(s,\alpha) d\alpha$ is a linear combination of $\zeta(s-1), \zeta(s-2),\ldots\ldots\ldots,\zeta(s-i)$

for each i=0,1,2,…………, N.

**Proof of corollary 8) :** We have $\int_0^1 \zeta(-m_1,\alpha)\cdot\zeta(-m_2,\alpha)\ldots\ldots\ldots\zeta(-m_r,\alpha)\cdot\zeta'(s,\alpha)d\alpha$

$= \int_0^1 P_N(\alpha)\cdot\zeta'(s,\alpha)d\alpha = \sum_{j=0}^{N} a_j \int_0^1 \alpha^j \cdot \zeta'(s,\alpha)d\alpha$ ,

where $P_N(\alpha)$ is the polynomial as in corollary 7) .

Note that $\int_0^1 \alpha^j \cdot \zeta'(s,\alpha)d\alpha = \int_0^1 \alpha^j d\left(\frac{\zeta'(s-1,\alpha)}{1-s} + \frac{\zeta(s-1,\alpha)}{(1-s)^2}\right)d\alpha$

$= \left[\alpha^j\left(\frac{\zeta'(s-1,\alpha)}{1-s} + \frac{\zeta(s-1,\alpha)}{(1-s)^2}\right)\right]_{\alpha=0}^1 - \int_0^1\left(\frac{\zeta'(s-1,\alpha)}{(1-s)} + \frac{\zeta(s-1,\alpha)}{(1-s)^2}\right)\cdot j\alpha^{j-1}d\alpha$

$= \frac{\zeta'(s-1)}{1-s} + \frac{\zeta(s-1)}{(1-s)^2} - j\int_0^1\left(\frac{\zeta'(s-1,\alpha)}{(1-s)} + \frac{\zeta(s-1,\alpha)}{(1-s)^2}\right)\cdot \alpha^{j-1}d\alpha$ .

We continue to integrate by parts j times .

Thus we find that $\int_0^1 \alpha^j \zeta'(s,\alpha)d\alpha$ is a linear combination of

$\zeta'(s-1), \zeta'(s-2), \ldots, \zeta'(s-j)$ ; $\zeta(s-1), \zeta(s-2), \ldots, \zeta(s-j)$

for each $j = 1,2,\ldots, N$ .

**Proof of Corollary 9)** : We have for Re s>1 ,

$\int_0^1 \zeta(0,\alpha)\zeta(1-s,\alpha)\cdot\zeta(2-s,\alpha)d\alpha = \frac{1}{2(s-1)}\int_0^1 (\frac{1}{2}-\alpha)\frac{\partial}{\partial\alpha}\zeta^2(1-s,\alpha)d\alpha$

$= \frac{1}{2(s-1)}\left\{(\frac{1}{2}-\alpha)\zeta^2(1-s,\alpha)|_{\alpha=0}^1 + \int_0^1 \zeta^2(1-s,\alpha)d\alpha\right\}$

$= \frac{1}{2(s-1)}\left(-\zeta^2(1-s) + \int_0^1 \zeta^2(1-s,\alpha)d\alpha\right) = \frac{1}{2(s-1)}\left(2(2\pi)^{-2s}\Gamma^2(s)\zeta(2s) - \zeta^2(1-s)\right)$ .

It is known that $\zeta(s,\alpha)$ is an analytic function of the complex variable s with a simple pole at s=1 . In author [1] , [2] it has been shown that $\zeta(s,\alpha)$ is an analytic function of the complex variable $\alpha$ except for possible singularities at non-positive integer values of $\alpha$ . For an integer $r \geq 0$, we write $\zeta^{(r)}(s,\alpha) = \dfrac{\partial^r}{\partial s^r}\zeta(s,\alpha)$. In author [1], [2] it has been shown that $\zeta^{(r)}(s,\alpha)$ is an analytic function of the complex variable $\alpha$ except for possible singularities at $\alpha = 0, -1, -2, \ldots\ldots$ .

We have $\zeta(s,\alpha) = \sum_{n\geq 0}(n+\alpha)^{-s}$ for Re s>1 ,

so that $(-1)^r \zeta^{(r)}(s,\alpha) = \sum_{n\geq 0}(n+\alpha)^{-s}\log^r(n+\alpha)$ for $r \geq 0$ .

Thus $(-1)^r \dfrac{\partial}{\partial \alpha}\zeta^{(r)}(s,\alpha) = r\sum_{n\geq 0}(n+\alpha)^{-s-1}\log^{r-1}(n+\alpha) - s\sum_{n\geq 0}(n+\alpha)^{-s-1}\log^r(n+\alpha)$.

This gives $\dfrac{\partial}{\partial \alpha}\zeta^{(r)}(s,\alpha) = -r\zeta^{(r-1)}(s+1,\alpha) - s\zeta^{(r)}(s+1,\alpha)$ for Re s>1 and $r \geq 1$ .

By analyticity of $\zeta^{(r)}(s,\alpha)$ as a function of complex variables $s$ and $\alpha$, the above result holds everywhere . This proves 2) of Proposition .

Let $k \geq 1$ be an integer . Let $\zeta_k(s,\alpha) = \sum_{n\geq k}(n+\alpha)^{-s}$ for Re s>1; and its analytic continuation . Then $\zeta_k(s,\alpha) = \zeta(s,\alpha) - \sum_{0\leq n\leq k-1}(n+\alpha)^{-s} = \zeta(s, k+\alpha)$ .

Let $\zeta_k(s) = \zeta(s) - \sum_{1\leq n\leq k-1}n^{-s}$ . Then in author [2], it has been shown that

$\zeta_k(s,\alpha) = \sum_{n\geq 0}\dfrac{(-\alpha)^n}{n!}s(s+1)\ldots\ldots(s+n-1)\cdot \zeta_k(s+n)$ in the disc $|\alpha| < k$,

where empty product stands for 1 .

This gives $\zeta(s,\alpha) = \sum_{0\leq n\leq k-1}(n+\alpha)^{-s} + \sum_{n\geq 0}s(s+1)\ldots\ldots(s+n-1)\zeta_k(s+n)\cdot \dfrac{(-\alpha)^n}{n!}$

$= \sum_{0\leq n\leq k-1}(n+\alpha)^{-s} + \phi_k(s,\alpha)$, say , where $\phi_k(s,\alpha)$ is an analytic function of α in the disc $|\alpha| < k$ with $k \geq 1$ arbitrary and $\phi_k^{(r)}(s,\alpha) = \dfrac{\partial^r}{\partial s^r}\phi_k(s,\alpha)$ can be obtained by term-by-term differentiation of $\phi_k(s,\alpha)$ with respect to s, r times.

Thus $\zeta^{(r)}(s,\alpha) = \sum_{0 \leq n \leq k-1}(n+\alpha)^{-s}(-\log(n+\alpha))^r + \phi_k^{(r)}(s,\alpha)$,

where $\phi_k^{(r)}(s,\alpha)$ is an analytic function of α in the disc $|\alpha| < k$ and $k \geq 1$ is arbitrary.

From the expression $\zeta(s,\alpha) = \sum_{0 \leq n \leq k-1}(n+\alpha)^{-s} + \phi_k(s,\alpha)$, where $\phi_k(s,\alpha)$ is analytic in the disc $|\alpha| < k$ with integer $k \geq 1$ arbitrary, we get the following :

1) If $s = -m$, where $m \geq 0$ is an integer, then $\zeta(-m,\alpha) = \sum_{0 \leq n \leq k-1}(n+\alpha)^m + \phi_k(-m,\alpha)$.

   Thus $\zeta(-m,\alpha)$ is an entire function of α.

2) If $s = m \geq 1$ is an integer, then $\zeta(m,\alpha) = \sum_{0 \leq n \leq k-1}(n+\alpha)^{-m} + \phi_k(m,\alpha)$.

   Thus $\zeta(m,\alpha)$ has a pole of order m at every non-positive integer value of α.

3) If Re s<0, then $\zeta^{(r)}(s,\alpha) = (-1)^r \sum_{0 \leq n \leq k-1}(n+\alpha)^{-s}\log^r(n+\alpha) + \phi_k^{(r)}(s,\alpha)$,

   where $\phi_k^{(r)}(s,\alpha)$ is analytic function of α in the disc $|\alpha| < k$. This shows that $\zeta^{(r)}(s,\alpha)$ is a continuous function of α, because if $\alpha = -n_0$, where $n_0$ is a fixed integer from amongst 0,1,2,……………, we have $\lim_{\alpha \to -n_0}(n_0+\alpha)^{-s}\log^r(n_0+\alpha) = 0$,

   for Re s<0. This also gives that if s is any complex number, the only possible singularities of $\zeta^{(r)}(s,\alpha)$ are non-positive integer values of α.

4) For Re s>1, we have $\zeta(s,\alpha) = \sum_{n \geq 0}(n+\alpha)^{-s}$.

   As $\dfrac{\partial}{\partial \alpha}\dfrac{\partial}{\partial s}(n+\alpha)^{-s} = \dfrac{\partial}{\partial s}\dfrac{\partial}{\partial \alpha}(n+\alpha)^{-s}$ for $\alpha \neq -n$,

   we have $\dfrac{\partial}{\partial \alpha}\dfrac{\partial}{\partial s}\left(\sum_{n \geq 0}(n+\alpha)^{-s}\right) = \dfrac{\partial}{\partial s}\dfrac{\partial}{\partial \alpha}\left(\sum_{n \geq 0}(n+\alpha)^{-s}\right)$ for Re s>1

and for $\alpha \neq -n$, where $n = 0,1,2,\ldots\ldots$ .

Thus $\dfrac{\partial}{\partial \alpha} \dfrac{\partial}{\partial s} \zeta(s,\alpha) = \dfrac{\partial}{\partial s} \dfrac{\partial}{\partial \alpha} \zeta(s,\alpha)$ for Re s>1 and for 0<α<1 .

Thus in view of analyticity of $\zeta(s,\alpha)$ as a function of each complex variable s and α , we have $\dfrac{\partial}{\partial \alpha} \dfrac{\partial}{\partial s} \zeta(s,\alpha) = \dfrac{\partial}{\partial s} \dfrac{\partial}{\partial \alpha} \zeta(s,\alpha)$ for $s \neq 1$ and $\alpha \neq 0,-1,-2,\ldots\ldots$ .

More generally , if $r_1, r_2 \geq 0$ are integers , we have

$$\dfrac{\partial^{r_1}}{\partial \alpha^{r_1}} \dfrac{\partial^{r_2}}{\partial s^{r_2}} \zeta(s,\alpha) = \dfrac{\partial^{r_2}}{\partial s^{r_2}} \dfrac{\partial^{r_1}}{\partial \alpha^{r_1}} \zeta(s,\alpha) \text{ for } s \neq 1 \text{ and } \alpha \neq 0,-1,-2,\ldots\ldots .$$

In particular , $\dfrac{\partial}{\partial \alpha}\left(\dfrac{\partial}{\partial s}\zeta(s,\alpha)\right)_{s=0} = \dfrac{\partial}{\partial s} \dfrac{\partial}{\partial \alpha} \zeta(s,\alpha)_{s=0}$ for $\alpha \neq 0,-1,-2,\ldots\ldots$ .

That is , $\dfrac{\partial}{\partial \alpha} \zeta'(0,\alpha) = \dfrac{\partial}{\partial s}(-s\zeta(s+1,\alpha))|_{s=0}$ .

Let $\zeta(s+1,\alpha) = \dfrac{1}{s} + \sum_{n \geq 0} \gamma_n(\alpha)s^n$ so that $s\zeta(s+1,\alpha) = 1 + \sum_{n \geq 0} \gamma_n(\alpha)s^{n+1} = \sum_{n \geq 0} \gamma_{n-1}(\alpha)s^n$ ,

where we have written $\gamma_{-1}(\alpha) = 1$

Thus $(s\zeta(s+1,\alpha))^{(r)}|_{s=0} = r!\gamma_{r-1}(\alpha)$ for $r = 0,1,2,\ldots\ldots$ .

Here the superscript (r) denotes r-th order derivative with respect to s .

Thus we have $\dfrac{\partial}{\partial \alpha} \zeta'(0,\alpha) = -\gamma_0(\alpha)$ .

However , we know $\zeta'(0,\alpha) = \log \dfrac{\Gamma(\alpha)}{\sqrt{2\pi}}$ for $\alpha \neq -n$, where $n = 0,1,2,\ldots\ldots$

Thus we have $\psi(\alpha) = -\gamma_0(\alpha)$ , where $\psi(\alpha) = \dfrac{\Gamma'}{\Gamma}(\alpha)$ .

Note that $\psi(1) = -\gamma_0(1) = -\gamma$ ,where $\gamma$ is Euler's constant .

The fact $(s\zeta(s+1,\alpha))^{(r)}|_{s=0} = r!\gamma_{r-1}(\alpha)$ gives $\left(\dfrac{\partial}{\partial\alpha}\zeta(s,\alpha)\right)^{(r)}|_{s=0} = -r!\gamma_{r-1}(\alpha)$.

Equivalently, we have $\dfrac{\partial}{\partial\alpha}\zeta^{(r)}(0,\alpha) = -r!\gamma_{r-1}(\alpha)$ for $r \geq 0$.

5) We have $\zeta'(0,\alpha) = -\sum_{0 \leq n \leq k-1}\log(n+\alpha) + \phi_k'(0,\alpha)$

so that $\dfrac{\partial}{\partial\alpha}\zeta'(0,\alpha) = -\sum_{0 \leq n \leq k-1}\dfrac{1}{n+\alpha} + \dfrac{\partial}{\partial\alpha}\phi_k'(0,\alpha)$.

That is $\psi(\alpha) = -\sum_{0 \leq n \leq k-1}\dfrac{1}{(n+\alpha)} + \dfrac{\partial}{\partial\alpha}\phi_k'(0,\alpha)$.

Thus $\psi(\alpha)$ has a simple pole at each non-positive integer value of α and consequently for $r \geq 1$, $\dfrac{\partial^r}{\partial\alpha^r}\psi(\alpha) = (-1)^{r-1}r!\zeta(r+1,\alpha).$ has a pole of order $(r+1)$ at each non-positive integer value of α and is analytic elsewhere.

This completes our proof.

### References

1] V.V. Rane, Instant Evaluation and demystification of $\zeta(n), L(n,\chi)$ that Euler, Ramanujan missed-I (ar Xiv org. website).

2] V.V. Rane, Instant Evaluation and demystification of $\zeta(n), L(n,\chi)$ that Euler, Ramanujan missed-II (ar Xiv org. website).

3] V.V.Rane, Functional Equation of Tornheim Zeta function and the α-Calculus of $\dfrac{\partial^r}{\partial s^r}\zeta(s,\alpha)$ (ar Xiv org. website).